\documentclass[a4paper]{amsart}
\usepackage{graphics}
\usepackage{amsmath}
\usepackage{latexsym}
\usepackage{amssymb}
\usepackage[all]{xy}
\xyoption{matrix}
\xyoption{arrow}

\begin{document}

\renewcommand{\th}{\operatorname{th}\nolimits}
\newcommand{\rej}{\operatorname{rej}\nolimits}
\newcommand{\extto}{\xrightarrow}
\renewcommand{\mod}{\operatorname{mod}\nolimits}
\newcommand{\Sub}{\operatorname{Sub}\nolimits}
\newcommand{\ind}{\operatorname{ind}\nolimits}
\newcommand{\Fac}{\operatorname{Fac}\nolimits}
\newcommand{\add}{\operatorname{add}\nolimits}
\newcommand{\Hom}{\operatorname{Hom}\nolimits}
\newcommand{\Rad}{\operatorname{Rad}\nolimits}
\newcommand{\RHom}{\operatorname{RHom}\nolimits}
\newcommand{\uHom}{\operatorname{\underline{Hom}}\nolimits}
\newcommand{\End}{\operatorname{End}\nolimits}
\renewcommand{\Im}{\operatorname{Im}\nolimits}
\newcommand{\Ker}{\operatorname{Ker}\nolimits}
\newcommand{\Coker}{\operatorname{Coker}\nolimits}
\newcommand{\Ext}{\operatorname{Ext}\nolimits}
\newcommand{\op}{{\operatorname{op}}}
\newcommand{\Ab}{\operatorname{Ab}\nolimits}
\newcommand{\id}{\operatorname{id}\nolimits}
\newcommand{\pd}{\operatorname{pd}\nolimits}
\newcommand{\A}{\operatorname{\mathcal A}\nolimits}
\newcommand{\C}{\operatorname{\mathcal C}\nolimits}
\newcommand{\D}{\operatorname{\mathcal D}\nolimits}
\newcommand{\X}{\operatorname{\mathcal X}\nolimits}
\newcommand{\Y}{\operatorname{\mathcal Y}\nolimits}
\newcommand{\F}{\operatorname{\mathcal F}\nolimits}
\newcommand{\Z}{\operatorname{\mathbb Z}\nolimits}
\renewcommand{\P}{\operatorname{\mathcal P}\nolimits}
\newcommand{\T}{\operatorname{\mathcal T}\nolimits}
\newcommand{\G}{\Gamma}
\renewcommand{\L}{\Lambda}
\newcommand{\bdot}{\scriptscriptstyle\bullet}
\renewcommand{\r}{\operatorname{\underline{r}}\nolimits}
\newtheorem{lem}{Lemma}[section]
\newtheorem{prop}[lem]{Proposition}
\newtheorem{cor}[lem]{Corollary}
\newtheorem{thm}[lem]{Theorem}
\newtheorem*{thmA}{Theorem}
\newtheorem*{thmB}{Theorem}
\newtheorem{rem}[lem]{Remark}
\newtheorem{defin}[lem]{Definition}


\title[Cluster-tilted algebras]{Cluster-tilted algebras of finite representation type}

\author[Buan]{Aslak Bakke Buan}
\address{Institutt for matematiske fag\\
Norges teknisk-naturvitenskapelige universitet\\
N-7491 Trondheim\\
Norway}
\email{aslakb@math.ntnu.no}

\author[Marsh]{Robert J. Marsh}
\address{Department of Mathematics \\
University of Leicester \\
University Road \\
Leicester LE1 7RH \\
England
}
\email{rjm25@mcs.le.ac.uk}

\author[Reiten]{Idun Reiten}
\address{Institutt for matematiske fag\\
Norges teknisk-naturvitenskapelige universitet\\
N-7491 Trondheim\\
Norway}
\email{idunr@math.ntnu.no}


\begin{abstract}
We investigate the cluster-tilted algebras of finite representation type over an algebraically
closed field. We
give an explicit description of the relations for the quivers for finite representation type.
As a consequence we show that a 
(basic) cluster-tilted algebra of finite type is uniquely determined by
its quiver. 
Also some necessary conditions on the shapes of quivers of cluster-tilted algebras of finite representation type
are obtained along the way.
\end{abstract}

\thanks{Aslak Bakke Buan was supported by a grant from the Norwegian Research Council.
}

\maketitle

\section*{Introduction}
Cluster categories $\C = \C_H$ associated with finite dimensional
hereditary algebras $H$ over a field $K$ (or more generally with $\Ext$-finite hereditary abelian 
$K$-categories with tilting object) were introduced in \cite{bmrrt}. 
An alternative description in Dynkin type $A$
was given in \cite{ccs1}.
The motivation came from
the theory of cluster algebras, introduced by Fomin and Zelevinsky in \cite{fz}, and
the connection between cluster algebras and quiver representations \cite{mrz}. A tilting
theory in cluster categories was developed in \cite{bmrrt}, and the associated endomorphism algebras 
$\End_{\C}(T)^{\op}$ for a (cluster-)tilting object $T$ in $\C$, were 
investigated in \cite{bmr1, bmr2}.
The cluster category has also motivated a Hall-algebra type definition of a cluster algebra
of finite type \cite{cc, ck}.

The cluster-tilted algebras are in spirit similar to the class of tilted algebras,
but their properties are quite different. On the one hand the tilted algebras
are nice from a homological point of view since they have global dimension at most two, 
while cluster-tilted algebras often have infinite global dimension. On the other hand, 
the indecomposable modules for a cluster-tilted algebra are in 1--1 correspondence 
with the indecomposable modules for the associated hereditary algebra. In particular,
a connected cluster-tilted algebra $\End_{\C_H}(T)^{\op}$
is of finite representation type if and only if $H$ is Morita equivalent to the path algebra
of a Dynkin quiver. Here we restrict to algebraically closed base fields $K$, and
thus to simply-laced Dynkin quivers.
Since we consider cluster-tilted algebras which are finite dimensional 
over an algebraically closed field, they
are (up to Morita equivalence) of the form $KQ/I$, where $Q$ is a finite quiver, and $I$
is some admissible ideal in the path algebra $KQ$, generated by a finite number of
paths.
The elements in $I$ are called {\em relations} if they are linear combinations 
$k_1 \rho_1 + \cdots + k_m \rho_m$ of paths $\rho_i$ in $Q$, all starting in the same vertex and
ending in the same vertex, and with each $k_i$ non-zero in $K$. 
If $m=1$, we call the relation a 
{\em zero-relation}. If $m=2$, we call it a commutativity-relation (and say that the paths $\rho_1$ and
$\rho_2$ commute).
For ease of notation we use the same symbol $\rho$ to denote a path, the corresponding element in the path algebra,
and the corresponding element in $KQ/I$.

A relation $\rho$ is called {\em minimal} if 
whenever $\rho = \sum_i \beta_i \circ \rho_i \circ \gamma_i$, where $\rho_i$ is a relation
for every $i$, then there is an index $i$ such that both $\beta_i$ and $\gamma_i$ are scalars.
When there is an arrow $i \to j$,
a path from $j$ to $i$ is called {\em shortest} if it contains no proper subpath
which is a cycle and if the full subquiver generated
by the induced oriented cycle contains no further arrows.
Our main result is the following.

\begin{thmA}
Let $Q$ be a finite quiver, $K$ an algebraically closed field and $I$ an ideal in the path algebra $KQ$, such that
$\G = KQ/I$ is a cluster-tilted algebra of finite representation type, and 
let $j$ and $i$ be vertices of $Q$.
\begin{itemize}
\item[(a)]{Assume there is an arrow $i \to j$. Then there are at most two shortest paths from
$j$ to $i$.}
\begin{itemize}
\item[(i)]{If there is exactly one, then this is a minimal zero-relation.}
\item[(ii)]{If there are two, $\rho$ and $\mu$, then $\rho$ and $\mu$ are not zero in $\G$,
they are disconnected, and there is a minimal relation $\rho + \lambda \mu$ for some 
$\lambda \neq 0$ in $K$.}
\end{itemize}
\item[(b)]{Up to multiplication by non-zero elements of $K$, there are no other minimal 
zero-relations or commutativity-relations.}
\item[(c)]{The ideal $I$ is generated by minimal zero-relations and minimal commutativity-relations.}
\end{itemize}
\end{thmA}

{\em Disconnected paths} are defined just after Lemma \ref{cycle-zero}.
Our main theorem has the following consequence.

\begin{thmB}
Let $\G$ be a cluster-tilted algebra of finite representation type over an algebraically closed field and 
with quiver $Q$. Let $I'$ be the ideal in $KQ$ with the 
following generators.
For $\alpha \colon i \to j$ choose a generator $\rho$ if there is exactly one shortest path $\rho$ from $j$ to $i$.
Choose a generator $\rho - \mu$ if there are two different shortest paths $\rho$ and $\mu$ from $j$ to $i$. 
Then $\G \simeq KQ/I'$.  
\end{thmB}


It is interesting to get a good understanding of the class of cluster-tilted algebras of finite 
representation type. In this paper we give some necessary conditions on the shape of the quiver for
such algebras, and use this to give an explicit description of the relations. 
In particular, it is worth noting that
our result means that the cluster-tilted algebras of finite type are (up to Morita equivalence) uniquely
determined by their quivers. It would be interesting to know to which extent this is true beyond finite
type. 
Our main theorem on describing the relations
of cluster-tilted algebras of finite type answers Conjecture 9.2 in \cite{bmrrt}. 
We remark that these relations appeared in \cite{ccs1}, and that
in \cite{ccs2} it is shown, independently, that they hold in a cluster-tilted algebra associated to
a (simply-laced) Dynkin quiver.
In type $A$ it is shown also that they are defining relations.
As an application of our main result, we can also complete the proof of Conjecture 1.1 from \cite{ccs1}.

The quivers occurring for cluster-tilted algebras associated with any hereditary finite dimensional algebra
are in 1--1 correspondence with the matrices occurring for an associated acyclic cluster algebra,
as shown in \cite{bmr2}. Hence obtaining results on the quivers of cluster-tilted algebras gives
information on matrices associated with cluster algebras. Note that
results on quivers of cluster algebras of finite type are given by A. Seven \cite{s}.
In fact Seven gives a list of all quivers which are cluster quivers of minimal infinite type; i.e. they 
define cluster algebras of 
infinite type,
but every full proper subquiver defines a cluster algebra of finite type.
For general notions in the representation theory of finite dimensional algebras we refer
to \cite{r} and \cite{ars}.

Some of the results in this paper have been presented at conferences in Uppsala (June 2004) and
Mexico (August 2004).

We would like to thank Claus Michael Ringel and \O yvind Solberg for some helpful conversations.

The first named author wishes to thank the Department of Mathematics at
the University of Leicester, and especially Robert J. Marsh, for their kind hospitality, during his stay there
in 2004. 

\section{Preliminaries}

\subsection{Cluster-tilted algebras and cluster categories}

Let $K$ be an algebraically closed field. We consider connected hereditary finite dimensional $K$-algebras.
Any such algebra $H$ is Morita equivalent to a path algebra $KQ$, for some finite
quiver $Q$. Furthermore, we assume $H = KQ$ is of finite representation type, that is, there 
is only a finite number of indecomposable objects, up to isomorphism, in the category $\mod H$ of 
finitely generated (left) $H$-modules. It is well known that this holds true if and only if
the underlying graph of $Q$ is a Dynkin graph. 

Let $\D = D^b(\mod H)$ be the bounded derived category. It is equipped with a shift-functor
$[1] \colon \D \to \D$. By \cite{hap} $\D$ has AR-triangles 
with corresponding translation-functor $\tau \colon \D \to \D$, with a quasi-inverse
$\tau^{-1}$. Let $F = \tau^{-1}[1]$ be the composition. It is an auto-equivalence on $\D$.
The cluster category is the orbit category $\C = \D/F$. The objects of $\C$ are the objects of
$\D$, while $\Hom_{\C}(A,B) = \amalg_i \Hom_{\D}(A,F^i B)$, see \cite{bmrrt}.

An object $T$ of $\C$ is called a (cluster-)tilting object if $\Ext^1_{\C}(T,T) = 0$ and $T$ is maximal
with respect to this property, i.e. if $\Ext^1_{\C}(T \amalg X,T \amalg X) = 0$, then
$X$ is a direct summand in a direct sum of copies of $T$.
Let $T = \bar{T} \amalg M$ be a basic tilting object, with $M$ indecomposable.
Then the following is shown in \cite{bmrrt}.

\begin{thm}\label{exchange}
\begin{itemize}
\item[(a)]{There is a unique indecomposable object $M^{\ast} \not \simeq M$ in $\C_H$ such that $\bar{T} 
\amalg M^{\ast}$ 
is a tilting object in $\C_H$.}
\item[(b)]{There are triangles $M^{\ast} \to B \to M \to $ and $M \to B' \to M^{\ast} \to $ in $\C_H$,
where $B \to M$ and $B' \to M^{\ast}$ are minimal right $\add \bar{T}$-approximations.}
\end{itemize}
\end{thm}

The endomorphism-algebra $\End_{\C}(T)^{\op}$ of a tilting object $T$ is called a {\em cluster-tilted algebra}.
The following was shown in \cite{bmr1}.

\begin{prop}\label{cta}
Let $\G = \End_{\C}(T)^{\op}$ be a cluster-tilted algebra with $\C = \C_H$ the cluster category
for some hereditary algebra $H$, and $T$ a tilting object in $\C$.
Then $\G$ is of finite representation type if and only if $H$ is of finite representation type.
In this case, the numbers of indecomposables in $\mod H$ and in $\mod \G$ are the same. 
\end{prop}

Assume $\G = \End_{\C}(T)^{\op}$ for a tilting object $T$ in $\C = \C_H$, where
$H$ is the path algebra of a Dynkin quiver $Q$. If the underlying graph
of $Q$ is the Dynkin graph $\Delta$, with $\Delta \in \{A_n, D_m, E_6, E_7, E_8 \}$
for $n \geq 1$ and $m \geq 4$, then we say that $\G$ is cluster-tilted of type $\Delta$.

A nice consequence of Proposition \ref{cta} is the following.

\begin{prop}\label{numbers}
If $\G = \End_{\C}(T)^{\op}$ is a connected cluster-tilted algebra of finite representation type,
then there is a unique Dynkin graph $\Delta$, 
such that $\G$ is cluster-tilted of type $\Delta$.
\end{prop}

\begin{proof}
Assume $\C = \C_H$ for some hereditary algebra $H$.
The number of indecomposable
objects for path algebras with underlying graph of type $A_n$ is $n(n+1)/2$.
For type $D_n$ (in case $n \geq 4$) it is $n(n-1)$, and for type $E_n$ (in case $6 \leq n \leq 8$) the numbers
are respectively $36,63$ and $120$. 
Combining this with Proposition \ref{cta}, and using that the number of simples for
$\G$ is the same as the number of simples of $H$, the claim follows.

\end{proof}

\subsection{Operations on quivers}

We consider two operations on quivers, factoring and mutating.
Given a quiver $Q$, we can remove a set $v_1, \dots, v_n$
of vertices, and all arrows starting or ending in any of the $v_i$.
The factor quiver is denoted $Q/\{v_1, \dots , v_n\}$.
Each vertex $v_i $ corresponds to a primitive idempotent $e_i$ in the path-algebra $\G =KQ$.
It is clear that 
$Q/\{v_1, \dots , v_n\}$ is the quiver of the algebra $\G/\G(e_1+ \cdots + e_n)\G$.
In \cite{bmr2} we proved the following.

\begin{prop}\label{factor}
Let $\G$ be a cluster-tilted algebra, and let $e$ be an idempotent. Then $\G/\G e \G$
is also a cluster-tilted algebra.
\end{prop}

We next consider {\em mutating}.
Let $Q$ be a quiver with no multiple arrows and with no loops and no oriented cycles of length two.
We describe mutation at a vertex $k$. 
The mutated quiver $Q'$ is obtained in the following way.
\begin{itemize}
\item[1.]{Add a new vertex $k^{\ast}$.}
\item[2.]{If there is a path $i \to k \to j$, then either:
\begin{itemize}
\item[I.]{If there is an arrow from $j$ to $i$, remove it.}
\item[II.]{If there is no arrow from $j$ to $i$, add an arrow from $i$ to $j$.}
\end{itemize} 
}
\item[3.]{For any vertex $i$ replace all arrows from $i$ to $k$ with arrows from $k^{\ast}$ to $i$,
and replace all arrows from $k$ to $i$ with arrows from $i$ to $k^{\ast}$.}
\item[4.]{Remove the vertex $k$.}
\end{itemize}

Note that in case $Q$ has two vertices, step 2 is void.
This definition can easily be extended to quivers with 
multiple arrows.
There is a canonical way to assign to a quiver with no loops and no oriented cycles of length two 
a square integral skew-symmetric matrix.
It can be easily seen that this definition is a special case of matrix mutation, as it appears in
the definition of cluster algebras \cite{fz}. In the case with no multiple arrows
we have integral skew-symmetric matrices where 
all elements are in $\{-1,0,1\}$. 

The following is a direct consequence of the main results in \cite{bmr2} and \cite{bmr1}.

\begin{prop}\label{mutating} 
Let $\G$ be a cluster-tilted algebra with quiver $Q$, and let $Q'$ be obtained from
$Q$ by a finite number of mutations.
Then there is a cluster-tilted algebra $\G'$ with quiver $Q'$.
Moreover $\G$ is of finite representation type if and only if $\G'$ 
is of finite representation type.
\end{prop}

Combining Propositions \ref{factor} and \ref{mutating}, we obtain the following.

\begin{prop}\label{combine}
Let $Q$ be the quiver of a cluster-tilted algebra, and let $Q'$ be
obtained from $Q$ by repeated mutating and/or factoring. Then $Q'$ is the quiver of a cluster-tilted algebra.
\end{prop}

The following also is a direct consequence of \cite{bmr2} and \cite{bmr1}.

\begin{lem}\label{maps}
Let $\G = KQ/I$ be cluster-tilted. Let $Q'$ be obtained from $Q$ by mutation, and let $\G'$ be the corresponding
cluster-tilted algebra. Let $Q''$ be obtained from $Q$ by factoring out a vertex $v$, and let $\G''$ be the corresponding
cluster-tilted algebra. 
\begin{itemize}
\item[(a)]{Let $\rho$ be a path in $Q$ such that $\rho$ also is a path in $Q'$, then it is a zero-path in $\G$ if and only
if it is a zero-path in $\G'$.}
\item[(b)]{Let $\rho$ be a path in $Q$ such that $\rho$ also is a path in $Q''$. If $\rho$ is a zero-path
in $\G$, then it also a zero-path in $\G''$. If $\rho$ is a non-zero path from $j$ to $i$ in $Q$, and there
is no path from $j$ to $i$ via $v$ in $Q$, then $\rho$ is non-zero in $\G''$.}
\item[(c)]{Assume there is a path $i \to k \to j$, and no arrow from $j$ to $i$. If we mutate at $k$, then the 
arrow $i \to j$ in $Q'$ represents the same map in the cluster category as the path $i \to k \to j$.}
\end{itemize}
\end{lem}

A crucial property for finite representation type is the following.

\begin{lem}\label{dimone}
Let $\G$ be a cluster-tilted algebra of finite representation type
and let $P_i, P_j$ be indecomposable projective
$\G$-modules. Then $\dim_k \Hom(P_i,P_j) \leq 1$.
\end{lem}

\begin{proof}
This follows from Lemma 8.2 in \cite{bmrrt}, using that any tilting object in a cluster category $\C$ is the
image of a tilting module for some algebra $H'$ with $\C_{H'} = \C$.
\end{proof}

\section{Double path avoiding quivers}\label{dpa-quivers}

In this section we show some necessary conditions on the quivers of cluster-tilted algebras of finite type.
It was shown in \cite{bmrrt,bmr2} that there are no loops and no (oriented) cycles of length at two in the quiver
of a cluster-tilted algebra.

A quiver $Q$ with no cycles of length two is called {\em double path avoiding} if 
\begin{itemize}
\item[-]{there are no multiple arrows in $Q$.}
\item[-]{the above property holds for any quiver $Q'$, obtained from $Q$ by possibly repeated factoring
and/or mutating.}
\end{itemize}

Any quiver of the form 
$$
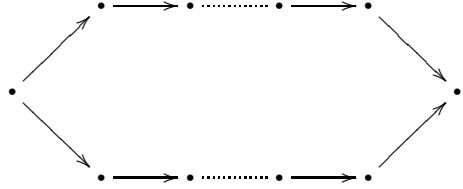
\xy
\xymatrix{
& \bdot \ar[r] & \bdot \ar@{.}[r] &  \bdot \ar[r] & \bdot \ar[dr] & \\
\bdot \ar[ur] \ar[dr] & & & & & \bdot \\
& \bdot \ar[r] & \bdot \ar@{.}[r] & \bdot \ar[r] & \bdot \ar[ur] &
}
\endxy
$$
is {\em not} double path avoiding.
This can be easily seen by picking any vertex $v$ except the source or the sink,
mutating at $v$, and then factoring out $v^{\ast}$. The path of length
two passing through $v$ is by this replaced by an arrow. 
$$
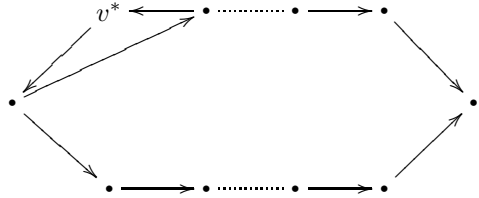
\xy
\xymatrix{
& v^{\ast} \ar[dl] & \bdot \ar[l] \ar@{.}[r] &  \bdot \ar[r] & \bdot \ar[dr] & \\
\bdot \ar[urr] \ar[dr] & & & & & \bdot \\
& \bdot \ar[r] & \bdot \ar@{.}[r] & \bdot \ar[r] & \bdot \ar[ur] &
}
\endxy
$$
We call this technique
{\em shortening of paths}. We can repeat this, removing all
vertices except the source and the sink, until we obtain the quiver
$$
\xy
\xymatrix{
\bdot  \ar@<1ex>[r] \ar@<-1ex>[r] & \bdot
}
\endxy
$$
More generally, we have the following.

\begin{lem}\label{oriented}
Let $Q$ be a quiver such that the underlying graph is a cycle of length $\geq 3$. Then $Q$ is double path
avoiding only if $Q$ is an oriented cycle.
\end{lem}

\begin{proof}
Assume that $Q$ is not an oriented cycle.
By applying shortening of paths, we can assume that the quiver
is alternating with $2n$ vertices, 
so all vertices are either sinks or sources.
$$
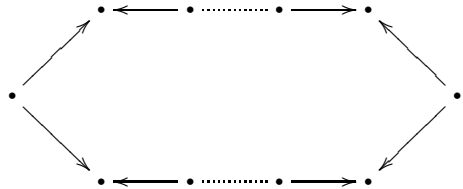
\xy
\xymatrix{
& \bdot & \bdot \ar[l] \ar@{.}[r] &  \bdot \ar[r] & \bdot & \\
\bdot \ar[ur] \ar[dr] & & & & & \bdot \ar[ul] \ar[dl] \\
& \bdot &  \bdot \ar[l] \ar@{.}[r] & \bdot \ar[r] & \bdot  &
}
\endxy
$$
Choose an arbitrary sink $v$, mutate at
$v$, and then apply shortening of paths on the new paths of length 2. If $n > 2$, then the resulting quiver
is alternating with $2n-2$ vertices. 
If $n=2$, then the reduced quiver is 
$$
\xy
\xymatrix{
\bdot \ar@<1ex>[r] \ar@<-1ex>[r] & \bdot
}
\endxy
$$ 
so by induction $Q$ is not double path avoiding.
\end{proof}

We note that this is shown in the finite type case in the cluster algebra context in
\cite{bgz}.
We will later see that if $Q$ is an oriented cycle of length $\geq 3$, then
$Q$ is double path avoiding. 

The motivation for the notion of double path avoiding quivers is the following.

\begin{prop}\label{finite-is-dpa}
Let $KQ/I$ be a cluster-tilted algebra of finite representation type.
Then $Q$ is double path avoiding.
\end{prop}

\begin{proof}
Assume $Q$ is not double path avoiding. Then there is a quiver $Q'$, obtained
by repeated factoring and mutating, such that $Q'$ has a double arrow.
Combining Propositions \ref{factor} and \ref{cta},
and using that there are no multiple arrows in the quiver for an algebra 
of finite type, we get a contradiction.
\end{proof}

Let us now define an important class of quivers.
For $a,b \geq 2$,
let $G(a,b)$ be the following quiver,
$$
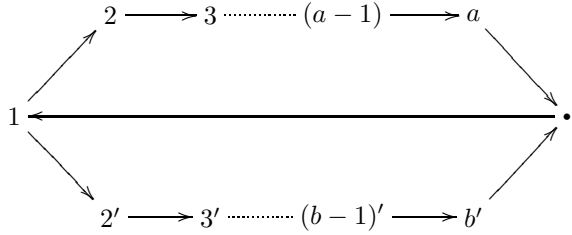
\xy
\xymatrix{
& 2 \ar[r] & 3 \ar@{.}[r] & (a-1) \ar[r] & a \ar[dr] & \\
1 \ar[dr] \ar[ur] & & & & & \bdot \ar[lllll] \\
& 2' \ar[r] & 3' \ar@{.}[r] & (b-1)' \ar[r] & b' \ar[ur] &
}
\endxy
$$
and let $T(a,b)$ be the quiver
$$
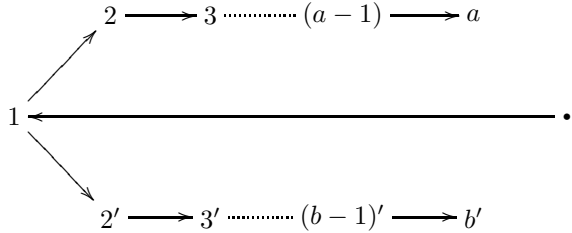
\xy
\xymatrix{
& 2 \ar[r] & 3 \ar@{.}[r] & (a-1) \ar[r] & a  & \\
1 \ar[ur]  \ar[dr] & & & & & \bdot \ar[lllll] \\
& 2' \ar[r] & 3' \ar@{.}[r] & (b-1)' \ar[r] & b'  &
}
\endxy
$$
\\
Two quivers are called {\em mutation-equivalent} if one can be obtained from
the other by a finite number of mutations. 

\begin{lem}\label{equiv}
The quivers $G(a,b)$ and $T(a,b)$ are mutation-equivalent (for $a,b \geq 2$).
\end{lem}

\begin{proof}
This can be seen by performing a series of mutations on $G(a,b)$, first mutate 
at $a$, and continue with $a-1$ and so on down to $2$. Next, mutate
at $b'$ and then $(b-1)'$, continue down to $2'$.
The resulting quiver is 
$$
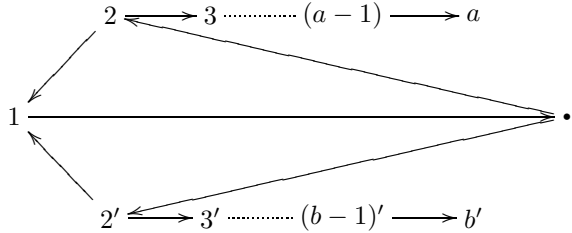
\xy
\xymatrix{
& 2 \ar[dl] \ar[r] & 3 \ar@{.}[r] & (a-1) \ar[r] & a  & \\
1 \ar[rrrrr]  & & & & & \bdot  \ar[llllu]  \ar[lllld]   \\
& 2'  \ar[ul]  \ar[r] & 3'  \ar@{.}[r] & (b-1)' \ar[r] & b' &
}
\endxy
$$
and it is clear that mutating at $1$ will give the desired
quiver $T(a,b)$.
\end{proof}

\begin{lem}\label{norel}
If an algebra with quiver $Q$ of Dynkin type is cluster-tilted,
then it is hereditary. 
\end{lem}

\begin{proof}
Assume $\G$ is cluster-tilted, with a quiver $Q$ of Dynkin type.
Assume there are relations.
They must all be zero-relations, that
is of the form $\alpha_1 \circ \alpha_2 \circ \cdots \circ \alpha_n$ for a path
$j_1 \overset{\alpha_1}{\rightarrow} j_2 \overset{\alpha_2}{\rightarrow} \cdots
j_{n-1} \overset{\alpha_n}{\rightarrow} j_n$.
We choose a zero-relation such that $n$ is minimal.
Factor out all vertices which are not on the path, to obtain the factor 
algebra $\G'$ with quiver $Q'$,
which is of type $A_n$, linear orientation.
$$
\xy
\xymatrix{
1 \ar[r] & 2 \ar[r] & 3 \ar[r] & \cdots \ar[r] & (n-1) \ar[r] & n
}
\endxy
$$
Then $\G'$ has a zero-relation. If $\G'$ was a cluster-tilted algebra of finite
representation type, then $\G'$ had to have the same number of indecomposable
modules as $A_n$, $D_n$ or $E_n$, by Proposition \ref{numbers}. 
But since $\G'$ is a proper factor algebra 
of $A_n$, this is impossible.  
\end{proof}

Then we have the following direct consequences of Proposition \ref{combine} and Lemma \ref{equiv}
and Lemma \ref{norel}.

\begin{thm}
Let $\G$ be a cluster-tilted algebra of finite representation type.
Assume the quiver of $\G$ is of type $G(a,b)$ with $a,b \geq 2$.
Then either $a = 2$, and $\G$ is cluster-tilted of type $D_{2+b}$, or
$(a,b)= (3,m)$ for $m \in \{3,4,5\}$, and  $\G$ is cluster-tilted of type $E_{3+m}$.
Moreover for these values of $(a,b)$, there is a cluster-tilted algebra with
quiver $G(a,b)$. 
\end{thm}

Let $C(n)$ be the cyclic quiver with $n$ vertices.
$$
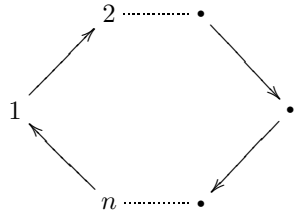
\xy
\xymatrix{
& 2 \ar@{.}[r] & \bdot \ar[dr] & \\
1 \ar[ur] & & & \bdot \ar[dl] \\
& n \ar[ul] \ar@{.}[r] & \bdot & }
\endxy
$$

\begin{prop}\label{cycles}
For any $n \geq 3$, there is a cluster-tilted algebra of type $D_n$ with quiver $C(n)$. 
(Here, we use the convention $A_3=D_3$.)
\end{prop}

\begin{proof}
For $n\geq 4$,
there is a cluster-tilted algebra with quiver
$G(2,n-2)$, that is
$$
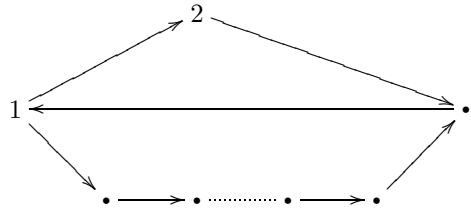
\xy
\xymatrix{
& & 2 \ar[drrr] & & & \\
1 \ar[urr] \ar[dr] & & & & & \bdot \ar[lllll] \\
& \bdot \ar[r] & \bdot \ar@{.}[r] & \bdot \ar[r] & \bdot \ar[ur] &
}
\endxy
$$
It is obtained by starting with $T(2,n-2)$, so it is of type $D_n$.
Now, mutate at $2$, to obtain $C(n)$.

For $n=3$, start with the quiver $A_3$, linear orientation
$$
\xy
\xymatrix{1 \ar[r] & 2 \ar[r] & 3
}
\endxy
$$
If we mutate at $2$, we obtain $C(3)$.
\end{proof}

Next we characterize the cluster-tilted algebras with cyclic quivers. Here
$\r$ denotes the Jacobson radical, the ideal generated by all arrows.

\begin{prop}\label{counting}
Let $\G$ be a cluster-tilted algebra with quiver $C(n)$.
Then $\G \simeq KC(n)/(\r)^{n-1}$.
\end{prop}

\begin{proof}

This is easy to check for $n=3$, since in this case there is (up to equivalence) only
one (connected) cluster category.

Now, let $n> 3$. Consider the composition $1 \to 2 \to \cdots \to j$, for $j \leq n-1$.
If this composition is zero, we can factor out the vertices $j+1, \cdots, n$, and we get 
a contradiction to Lemma \ref{norel}.
By Propositions \ref{numbers} and \ref{cycles}, 
the cluster-tilted algebra $\G$ must have the same number of indecomposables
as $D_n$. It is well known that this number is $n(n-1)$. 
Now we use the fact that $\G$ is a Nakayama algebra, so the number of 
indecomposable modules can be easily computed (see e.g. \cite{ars}, chapter IV.2). It is the sum of
the dimensions of the indecomposable projectives. Since any composition of $n-2$ arrows is nonzero by
Proposition \ref{norel},
each projective has dimension at least $n-1$. There are $n$ indecomposable projectives,
so it is clear that in fact the length of each of them must be exactly $n-1$, since the sum is $n(n-1)$.
\end{proof}

\begin{prop}\label{commutativity}
Given a cluster-tilted algebra with quiver $G(a,b)$.
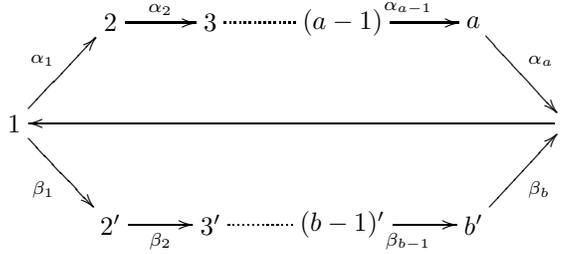
$$
\xy
\xymatrix{
& 2 \ar[r]^{\alpha_2} & 3 \ar@{.}[r] & (a-1) \ar[r]^{\alpha_{a-1}} & a \ar[dr]^{\alpha_a} & \\
1 \ar[dr]_{\beta_1} \ar[ur]^{\alpha_1} & & & & & \bdot \ar[lllll] \\
& 2' \ar[r]_{\beta_2} & 3' \ar@{.}[r] & (b-1)' \ar[r]_{\beta_{b-1}} & b' \ar[ur]_{\beta_b} &
}
\endxy
$$
Let $\alpha = \alpha_1 \circ \alpha_2 \circ \cdots \circ \alpha_a$ and
$\beta = \beta_1 \circ \beta_2 \circ \cdots \circ \beta_b$.
Then $\alpha \neq 0$ and $\beta \neq 0$ considered as elements in the cluster-tilted algebra.
Furthermore, there is a non-zero scalar $k$ such that
$\alpha + k \beta = 0$ in the cluster-tilted algebra.
\end{prop}

\begin{proof}
By applying shortening of paths, and using Lemma \ref{maps}, we may assume that $a=b=2$, so we have the quiver
$$
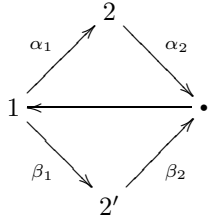
\xy
\xymatrix{
& 2 \ar[dr]^{\alpha_2} & \\
1 \ar[dr]_{\beta_1} \ar[ur]^{\alpha_1} & & \bdot \ar[ll] \\
& 2' \ar[ur]_{\beta_2} &
}
\endxy
$$
If we mutate at the vertex $2$, we obtain the cycle $C(4)$,
thus by Proposition \ref{counting} $\beta =  \beta_1 \circ \beta_2 \neq 0$. By symmetry we
have also $\alpha \neq 0$. 

If we instead of mutating at $2$, do factoring
at $2$, we get the quiver $C(3)$. It is therefore clear that
in the corresponding factor algebra the image of $\beta$ is zero. 
Thus, the assertion of the proposition follows.
\end{proof}

A subquiver $Q'$ of $Q$ is called full if every arrow
$i \to j$ in $Q$, with $i,j$ in $Q'$, is also an arrow in $Q'$.  
The subquiver generated by a collection of paths $\alpha_1, \dots, \alpha_m$ in 
$Q$ is the smallest
full subquiver of $Q$ containing the paths.
We say that an oriented cycle is {\em pure} if no proper subpath is a cycle.
We say that a path is full if the subquiver generated by the path
contains only the arrows on the path.
Assume there is an arrow $i \to j$. Then a path from $j$ to $i$ is called {\em shortest} 
if the induced oriented cycle is full and pure.

Propositions \ref{counting} and \ref{commutativity} have the following consequence,
which gives an alternative proof of Proposition 3.5 in \cite{ccs2}.

\begin{cor}\label{shortest}
Assume there is an arrow $i \to j$.
If there is exactly one shortest path $\alpha$
from $j$ to $i$,
then $\alpha = 0$. 
If there are exactly two shortest paths $\alpha$ and $\beta$ from $j$ to $i$,
then $\alpha + k \beta= 0$ for some scalar $k \neq 0$.
\end{cor}

\begin{proof}
The first part follows directly from Proposition \ref{counting}.
If there are exactly two shortest paths, the full subquiver generated by these two paths
must be as in Proposition \ref{commutativity}.
\end{proof}

Lemma \ref{dimone}
has the following direct consequences.

\begin{lem}\label{no-forward}
Assume there is an arrow $i \to j$ in the quiver of a cluster-tilted algebra
$\G$ of finite representation type. If $\beta$ is a path of length $>1$ from $i$ to $j$,
then $\beta = 0$. 
\end{lem}

\begin{lem}\label{cycle-zero}
Let $\alpha$ be an oriented cycle in the quiver of a cluster-tilted algebra
$\G$ of finite representation type. Then $\alpha = 0$.
\end{lem}
 
We say that two paths from $j$ to $i$ are {\em disjoint} if they
have no common vertices except $j$ and $i$. 
We say that two disjoint paths 
$\alpha = \alpha_1 \circ \alpha_2 \circ \cdots \circ \alpha_a$ and
$\beta = \beta_1 \circ \beta_2 \circ \cdots \circ \beta_b$ 
are {\em disconnected} if the full subquiver generated by
the paths contains no further arrows except possibly an arrow from $i$ to $j$.

For the rest of this section we assume that $\G = KQ/I$ is a cluster-tilted 
algebra of finite representation type, and thus
$Q$ is a quiver which is double path
avoiding.
There are several properties needed for the proof of the main result.
We state these as a series of lemmas.

\begin{lem}\label{alternating}
\begin{itemize}
\item[(a)]{Assume there is subquiver $Q'$ of $Q$ which is a non-oriented cycle. Let $v$ be a sink (source) of $Q'$.  
Then in the full subquiver generated by $Q'$, there is an arrow starting (ending) in $v$.}
\item[(b)]{Any arrow lying on a (not necessarily oriented) cycle, also lies on a full oriented cycle.}
\end{itemize}
\end{lem}

\begin{proof}
Part (a) follows straightforwardly from Lemma \ref{oriented}. For (b) one
can use the same lemma and induction on the length of the cycle.
\end{proof}

\begin{lem}\label{full-cycles}
Assume there is an oriented cycle $$j_1 \overset{\alpha_1}{\rightarrow} 
j_2 \overset{\alpha_2}{\rightarrow} \cdots \rightarrow j_{n-1} \overset{\alpha_{n-1}}{\rightarrow} j_n 
\overset{\alpha_n}{\rightarrow} j_1$$
in $Q$ such that the composition $\alpha_1 \circ \alpha_2 \circ \cdots \circ \alpha_{n-1} \neq 0$.
Then the cycle is full and pure.
\end{lem}

\begin{proof}
Note that by Lemma \ref{cycle-zero}, the cycle is pure.
Consider the cycle
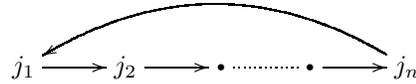
\vspace{1cm}
$$
\xy
\xymatrix{j_1 \ar[r] & j_2 \ar[r] & \bdot \ar@{.}[r] & \bdot \ar[r] & j_n \ar@/_2pc/[llll] \\   
}
\endxy
$$
\\
Let $Q'$ be the full subquiver generated by this cycle.
By Lemma \ref{no-forward}, there are no arrows $j_a \to j_b$ for $b>a+1$ in $Q'$.
Assume there is an arrow $j_b \to j_a$ for $a<b-1$, and choose such an arrow with
$a$ minimal. Consider the factor quiver $Q''$, obtained by deleting $j_{a+1}, \dots, j_{b-1}$. 
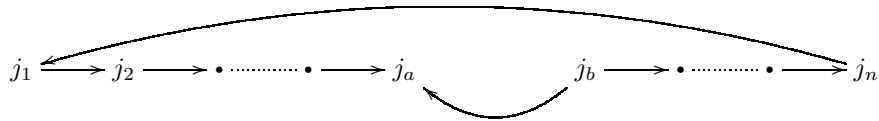
\vspace{1cm}
$$
\xy
\xymatrix{j_1 \ar[r] & j_2 \ar[r] & \bdot \ar@{.}[r] & \bdot \ar[r] & j_a  &   & j_b \ar[r] \ar@/^1.5pc/[ll] & 
  \bdot \ar@{.}[r] & \bdot \ar[r] & j_n \ar@/_2pc/[lllllllll] \\
}
\endxy
$$
\\
\\
By Lemma \ref{alternating}, there must be an arrow in $Q''$ starting in $j_a$. 
By choice of $a$, this arrow must end in some $j_{b'}$ for $b' > b$. This is a contradiction to
Lemma \ref{no-forward}.
\end{proof}

\begin{lem}\label{two-dis}
Assume there are two disjoint non-zero paths $\rho, \mu$ from $j$ to $i$, and an arrow
$i \to j$. Then the subquiver generated by $\rho$ and $\mu$  
contains no further arrows.
\end{lem}

\begin{proof}
Consider the quiver 

$$
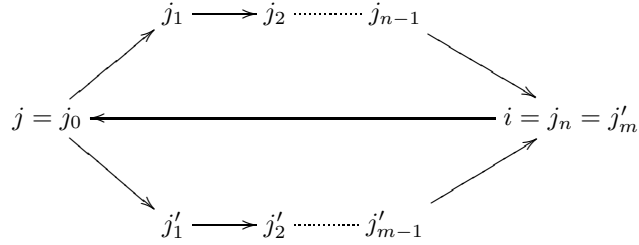
\xy
\xymatrix{
& j_1 \ar[r] & j_2 \ar@{.}[r] & j_{n-1} \ar[dr] & \\
j= j_0  \ar[dr] \ar[ur] & & & & i = j_n= j'_m \ar[llll] \\
& j'_1 \ar[r] & j'_2 \ar@{.}[r] & j'_{m-1} \ar[ur] &
}
\endxy
$$
where the upper path from $j$ to $i$ is $\rho$, the lower is $\mu$. 
If there is an additional arrow, then
by Lemma \ref{full-cycles}, it must start in $j_a$ for some $0<a <n$, and end
in $j'_b$ for some $0<b<m$ (or opposite). Such arrows are called 
crossing arrows.
If a crossing arrow exists, there is one $j_a \to j'_b$ such that the subquiver $Q'$ generated
by the subpath $j_0 \to \cdots \to j_a$ of $\rho$, the subpath $j_0 \to \cdots \to j'_b$ of $\mu$, and the
arrow $j_a \to j'_b$ contains no further crossing arrows. In $Q'$, there must, by
Lemma \ref{alternating}, be an arrow starting in $j'_b$. But by Lemma \ref{full-cycles}, only
crossing arrows starting in $j'_b$ can exist, so we have a contradiction.
\end{proof}

\begin{lem}\label{two-paths}
Assume there are two disconnected paths $\alpha, \beta$  from $j$ to $i$, and assume both are non-zero.
Then there is an arrow from $i$ to $j$.
\end{lem}

\begin{proof}

This follows directly from Lemma \ref{alternating}.

\end{proof}

%

\section{Minimal relations}

In this section we analyse minimal relations in cluster-tilted algebras
and the corresponding paths in their quivers. This will be used in the 
next section to prove our main result.
The following useful result of Assem-Br{\"u}stle-Schiffler and Reiten-Todorov gives
a major simplification of our original proof in this section. 
We include a proof for the convenience of the reader.
See also \cite{kr} for related considerations.

\begin{lem}\label{rt-abs}
Let $\G$ be a cluster-tilted algebra over $k$ and
let $S,S'$ be two simple $\G$-modules. Then $\dim_K \Ext^1(S',S) \geq \dim_K \Ext^2(S,S')$.
\end{lem}

\begin{proof}
We have $\G = \End_{\C}(T)^{\op}$ where $T$ is a tilting object in a cluster category $\C = \C_H$. Write
$T = \bar{T} \amalg M$, so that $\Hom_{\C}(T,M)$ is the projective cover of
$S$, and choose an indecomposable $M^{\ast} \not \simeq M$ such that $\bar{T} \amalg M^{\ast}$
is a tilting object. This is possible according to Theorem \ref{exchange}, and in the notation
of this theorem we have exchange triangles $M^{\ast} \to B \to M \to$ 
and $M \to B' \to M^{\ast} \to$, and hence exact sequences $$(T,M^{\ast}) \to (T,B) \to (T,M) \to S \to 0$$
and $$(T,M) \to (T,B') \to (T,M^{\ast}) \to (T,M[1]) = 0.$$ 
This gives the start of a projective resolution $$(T,B') \to (T,B) \to (T,M) \to S \to 0$$ of the simple
$\G$-module $S$, and the start of a minimal projective resolution $$(T,B'_0) \to (T,B) \to (T,M) \to S \to 0,$$
where $B'_0 \mid B'$. Let $(T,M') \to S'$ be a projective cover. Then $\dim_K \Ext^2(S,S') = \dim_K
\Hom((T,B'_0),S')$ is the multiplicity of $M'$ in $B'_0$.

On the other hand, since $M \to B'$ is a minimal left $\add \bar{T}$-approximation,
we have that $\dim_K \Ext^1(S',S)$
is the multiplicity of $M'$ in $B'$. This gives the desired inequality.
\end{proof}

Recall that a relation $\rho$ is called minimal if 
whenever $\rho = \sum_i \beta_i \circ \rho_i \circ \gamma_i$, where $\rho_i$ is a relation
for every $i$, then there is an index $i$ such that both $\beta_i$ and $\gamma_i$ are scalars.
We then have the following.

\begin{prop}\label{arrow-back}
Let $\rho$ be a path from $j$ to $i$
in $Q$ such that $\rho$ is a minimal zero-relation. Then there is an arrow $i \to j$, and
the induced oriented cycle is full and pure.
\end{prop}

\begin{proof}
Since there is a minimal zero-relation $\rho$ given by a path from $j$ to $i$, we know by \cite{Bong}
that $\Ext^2(S_j,S_i) \neq 0$, and hence $\Ext^1(S_i,S_j) \neq 0$ by Lemma \ref{rt-abs}. Hence there is an arrow from $i$
to $j$. We now want to show that the cycle induced by the minimal zero-relation $\rho$ from $j$ to $i$ and 
the arrow
from $i$ to $j$ is pure. 
Assume $\rho$ is given by the path $j= j_0 \to j_1 \to \cdots \to j_n = i$.
In case $(a,b) \neq (0,n)$ it follows from
Lemma \ref{no-forward} that there are no arrows $j_a \to j_b$ for $b> a+1$.
There are no arrows $j_0 \to j_n$ since there are no oriented cycles of length two.
Assume there was an arrow in the opposite direction. Choose such an arrow $\alpha \colon j_b \to j_a$ such that
$b$ is maximal.
Then consider the factor quiver obtained by factoring out $j_{a+1}, \cdots, j_{b-1}$, and all
vertices not on $\rho$. 
In this quiver $j_b$ is a source, since a potential arrow ending in $j_b$ would be an arrow of
the form $j_a \to j_b$ for $b > a+1$. Thus the underlying graph of the quiver contains a non-oriented
cycle, and by Lemma \ref{alternating} we have a contradiction.
Hence the cycle is full. Since $\rho$ is a minimal zero-relation it is clear that the cycle is pure
by Lemma \ref{cycle-zero}.
\end{proof}

\begin{lem}\label{maxtwo}
If there is an arrow $i \to j$, then there are at most two non-zero paths from $j$ to $i$.
If there are exactly two paths, they are disconnected. 
\end{lem}

\begin{proof}
Let $\rho \colon j=j_0 \to j_1 \to j_2 \to \dots \to j_n=i$ 
be a non-zero path. Assume there is a second non-zero path 
$\mu \colon j=j_0' \to j_1' \to j_2' \to \dots \to j_m'=i$. We show that
$\rho$ and $\mu$ are disconnected. 
Assume first that $\rho$ and $\mu$ have a common vertex $j_t = j'_{t'}$ for $1<t<n$ and $1<t'<m$.
Let $l>t$ be minimal such that $j_l$ is also on $\mu$.
Consider the subquiver $Q'$ generated by the two paths from 
$j_t$ to $j_l$.
Since $Q'$ is double path avoiding, there must be an arrow in $Q'$ starting in 
$j_l$, by Lemma \ref{alternating}. This contradicts Lemma \ref{full-cycles}.
\\
\\
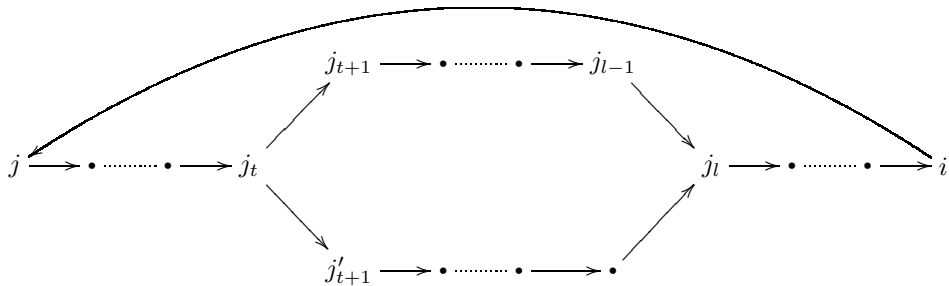
$$
\xy
\xymatrix@C=19pt{
& & & &  j_{t+1 }\ar[r] & \bdot \ar@{.}[r] & \bdot \ar[r] & j_{l-1} \ar[dr] & & & & \\
j \ar[r] & \bdot \ar@{.}[r] & \bdot \ar[r] & j_t \ar[ur] \ar[dr] & & & & & j_l \ar[r] & \bdot \ar@{.}[r] 
& \bdot \ar[r] & i \ar@/_5pc/[lllllllllll] \\
& & & & j'_{t+1} \ar[r] & \bdot \ar@{.}[r] & \bdot \ar[r] & \bdot \ar[ur] & & & & 
}
\endxy
$$
\\
So the two paths are disjoint. 
It now follows from Lemma \ref{two-dis} that they are also disconnected.

Now assume that there are three pairwise disconnected paths.
We first consider the quiver $Q$
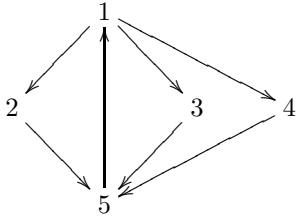
$$
\xy
\xymatrix{
& 1 \ar[dl] \ar[dr] \ar[drr] & & \\
2 \ar[dr] & & 3 \ar[dl] & 4 \ar[dll] \\
& 5 \ar[uu] & & 
}
\endxy
$$
If we mutate $Q$ at the vertex $4$, and then factor out $4^{\ast}$,
we obtain the following quiver
$$
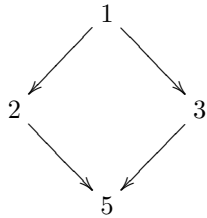
\xy
\xymatrix{
& 1 \ar[dl] \ar[dr] &  \\
2 \ar[dr] & & 3 \ar[dl]  \\
& 5  & & 
}
\endxy
$$
Hence $Q$ is not double path avoiding.
If $Q$ has paths of length $>2$, we can use shortening of paths to
reduce to the above case.
So there are at most two disconnected cycle-free paths, and thus at most two non-zero paths. 
\end{proof}

It follows from Lemma \ref{dimone} that when $\rho$ and $\mu$ are non-zero paths from $j$ to $i$,
then there is a commutativity-relation $\rho + k \mu$ for some non-zero scalar $k$.
Note that a relation is minimal in our definition if and only if
it is not in $\r I + I \r$, where
$\r$ denotes the Jacobson radical of $kQ$. 
We first characterize the minimal commutativity-relations.

\begin{lem}\label{only-disc}
Let $\rho$ and $\mu$ be distinct non-zero paths from $j$ to $i$.
Then the relation $\rho + k \mu$ with $k$ a non-zero scalar is minimal  
if and only if $\rho$ and $\mu$ are disconnected.
\end{lem}

\begin{proof}
Assume first there is a relation $\rho + k \mu$, with $\rho$ and $\mu$ 
both non-zero, which is not disconnected. We want to show that it is not minimal.

We first assume that the paths are disjoint.
There are then two cases, either any additional arrow is a back-arrow,
or there are crossing arrows.
We first assume that there are no crossing arrows, so the only additional arrows are
back-arrows on either $\rho$ or $\mu$. If there is an arrow $i \to j$, there are no further back-arrows by
Lemma \ref{two-dis}, so in this case the paths are disconnected. Consider the path $\rho$. Since 
it becomes zero if we factor out one of the vertices in $\mu$, it is clear that the path is not full.
That is, there are additional arrows on $\rho$, which must be back-arrows. By symmetry there are additional back-arrows
on $\mu$.
Since $Q$ is double path avoiding there must by Lemma \ref{alternating} be a back-arrow starting in $i$. 
Assume it ends in a vertex $j_a$
on $\rho$. 
Now consider $\mu$. There must also be a back-arrow starting in some $j'_{b'}$ on $\mu$. By removing 
vertices, we can 
obtain a subquiver which is
a non-oriented cycle, with $j'_{b'}$ a source, thus there must be an arrow ending in 
$j'_{b'}$ by Lemma \ref{alternating}.
This arrow must start in $j'_{c'}$ for some $c' > b'$.  
Continuing, we get a path of back-arrows on $\mu$ from $i$ to $j'_{b'}$.
Assume the first arrow on this path is $i \to j'_{a'}$. 
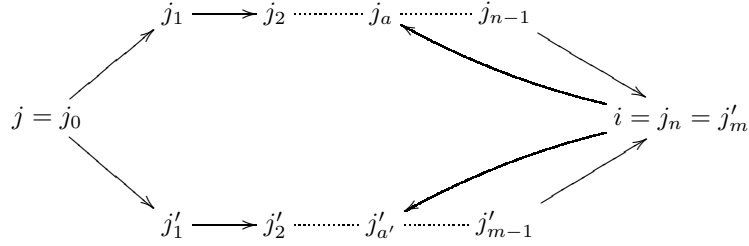
$$
\xy
\xymatrix{
& j_1 \ar[r] & j_2 \ar@{.}[r] & j_a \ar@{.}[r] & j_{n-1} \ar[dr] & \\
j= j_0  \ar[dr] \ar[ur] & & & & & i = j_n= j'_m \ar@/^0.5pc/[ull] \ar@/_0.5pc/[dll] \\
& j'_1 \ar[r] & j'_2 \ar@{.}[r] & j'_{a'} \ar@{.}[r] & j'_{m-1} \ar[ur] &
}
\endxy
$$
\\
Consider now the subquiver obtained by deleting the vertices between $j_a$ and $i$ on $\rho$ and between $j'_{a'}$ and $i$
on $\mu$. By Lemma \ref{alternating} there must be an arrow in this subquiver which ends in $i$. This arrow must
be a short-arrow on either $\rho$ or $\mu$, which gives a contradiction.

Next we keep the assumption that the paths are disjoint, but now assume crossing arrows.
Then there is a crossing arrow $j_t \to j'_{s'}$ such that the there are no crossing arrows between 
$\delta \colon j=j_0 \to j_1 \to \dots \to j_t \to j'_{s'}$ and the path 
$\mu_{s} \colon j=j_0 \to j'_1 \to \dots \to j'_{s'-1} \to j'_{s'}$ (or there is an arrow from a vertex
on $\mu$ to a vertex on $\rho$ with this property). 
$$
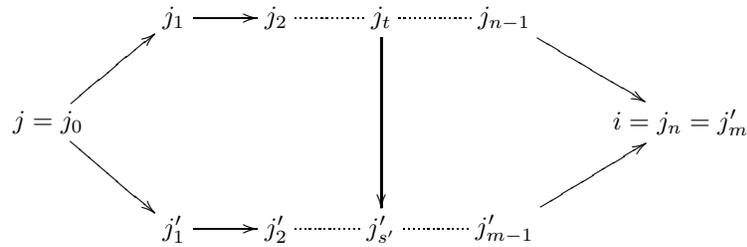
\xy
\xymatrix{
& j_1 \ar[r] & j_2 \ar@{.}[r] & j_t \ar[dd] \ar@{.}[r] & j_{n-1} \ar[dr] & \\
j= j_0  \ar[dr] \ar[ur] & & & & & i = j_n= j'_m \\
& j'_1 \ar[r] & j'_2 \ar@{.}[r] & j'_{s'} \ar@{.}[r] & j'_{m-1} \ar[ur] &
}
\endxy
$$
\\
We claim that $\delta$ is non-zero. Assume to the contrary that $\delta = 0$.
Then there must be back-arrows on $\delta$ by Lemma \ref{norel}.
Especially one can show that there must be a back-arrow on $\delta$ starting in $j'_{s'}$.
This can be seen by starting with a back-arrow starting in some $j_x$, and then noting that there
must also be a back-arrow ending in $j_x$. Continuing as for $\mu$ above, we obtain the required back-arrow.

The back-arrow on $\delta$ starting in $j'_{s'}$ must end in $j$, by the choice of the crossing arrow $j_t \to j'_{s'}$.
It is now straightforward to see that there can be no other back-arrows on $\delta$.
By Lemma \ref{full-cycles} there are no further back-arrows on $\mu_s$, in addition
to the arrow $j'_{s'} \to j$. Hence by Lemma \ref{commutativity} the path $\delta$ is non-zero, a
contradiction.

Hence $\delta$ is non-zero, and thus commutes with $\mu_s$.

Consider the subpath $\mu_e \colon j'_{s'} \to \cdots \to i$ of $\mu$ and let
$\delta'$ be the composition of $j_t \to  j'_{s'}$ with $\mu_e$.
Consider the subpath $\rho' \colon j_t \to \cdots \to i$ of $\rho$.
Then $\delta \circ \mu_e = k \mu$ for some non-zero $k$ and thus  $\delta \circ \mu_e$ is non-zero. 
Especially $\delta'$ is non-zero and thus commutes
with $\rho'$.
This shows that in this case the relation $\rho + k \mu$ is not minimal.

Now assume $\rho$ and $\mu$ are not disjoint.
There is some $t >0$, such that $j_t$ is on both $\rho$ and $\mu$.
The subpaths $\rho' \colon j \to \cdots \to j_t$ and $\mu' \colon j \to \cdots \to j_t$ of $\rho$ and $\mu$ respectively,
commute. The same is true for the subpaths starting in $j_t$. Thus $\rho + k \mu$ is not minimal also
in this case.

Now, assume $\rho, \mu$ are
non-zero disconnected paths from $j$ to $i$. 
We need to show that the relation $\rho +k \mu$ is minimal. 
Assume  $\rho +k \mu = \sum_i \beta_i \circ \tilde{\rho}_i \circ \gamma_i$,
where the $\tilde{\rho_i}$ are relations. 
One of the terms on the right hand side, say  $\beta_s \circ \tilde{\rho}_s \circ \gamma_s$,
must have $\rho$ as a summand. Assume $\rho = \beta_s \circ \tilde{\rho}_s^1 \circ \gamma_s$,
where $\tilde{\rho}_s^1$ is a subpath of $\rho$. Since $\tilde{\rho}_s^1$ is involved in a relation, it is zero in some
factor algebra. Since $\rho \neq 0$ there is no short-arrow, 
and hence there must be a back-arrow on $\tilde{\rho}_s^1$. But, by the assumption on $\rho$ this means
that $\tilde{\rho}_s^1 = \rho$. Hence $\beta_s$ and $\gamma_s$ are scalars, and we are done.
\end{proof}

We want to show that the minimal commutativity-relations together with 
the minimal zero-relations generate $I$. 
We have the following preliminary result.

\begin{lem}\label{sum}
Any minimal relation $\rho$ is of the form $\rho_m + \rho_c$, where $\rho_m$
is either a minimal zero-relation or a minimal commutativity-relation and $\rho_c$ is a relation that is not minimal.
\end{lem}

\begin{proof}
Assume we have a linear combination $\rho = k_1 \rho_1 + \cdots + k_n \rho_n$ of distinct paths
between $j$ and $i$ in $KQ$ which is a minimal relation. We clearly can assume that none of the $\rho_i$ are non-minimal
zero-paths, because they can be absorbed into $\rho_c$.

First assume that $\rho_1$ is a zero-path. Then, since it is minimal, there is an arrow 
$i \to j$, by Proposition \ref{arrow-back}. 
Note now that any path $\rho_i$ is full and gives rise to a pure cycle.
For $\rho_i \neq  0$ this follows from Lemma \ref{full-cycles}, otherwise
it follows from Proposition \ref{arrow-back}.

Assume now there are at least two distinct 
paths $\mu \colon j= j_0 \to j_1 \to \cdots \to j_n = i$ and 
$\mu' \colon j=j'_0 \to j'_1 \to \cdots \to j'_{n'} = i$ from $j$ to $i$, 
each giving rise to a full and pure cycle.
We first show that $\mu$ and $\mu'$ are disjoint. If not we have
a pair of vertices $(s,t) \neq (j,i)$ on $\mu$ and $\mu'$ and two disjoint paths $\mu_p$ and $\mu'_p$
from $s$ to $t$.
Assume $s= j_u = j'_{u'}$ and $t= j_v = j'_{v'}$. 
By assumption, there are no back-arrows on $\mu_p$ or $\mu'_p$.
We show that there are no crossing arrows. Assume there are crossing arrows, and let
$w$ with $u<w<v$ be minimal such that there is a crossing arrow starting or ending in $j_w$.
Then choose the smallest $w'$ with
$u' < w' < v'$ such that there is a crossing arrow $\beta$ between $j_w$ and $j'_{w'}$.
Consider the case where $\beta$ starts in $j_w$ and ends in $j'_{w'}$.
Let $\delta$ denote the induced path from $j_u$ to $j'_{w'}$ via $j_w$. 
Consider the subquiver generated by the two paths from 
$j_u$ to $j'_{w'}$. By choice of crossing arrow, these paths are disconnected. 
Since we have a cycle with an alternating vertex in $j'_{w'}$, this is a contradiction.
The case where $\beta$ starts in $j'_{w'}$ and ends in $j_{w}$ is similar. 
Hence, the paths $\mu_p$ and $\mu'_p$ are disconnected. Thus, by Lemma \ref{alternating}
there must be an arrow from $t$ to $s$,
but this contradicts the assumption that $\mu$ (and $\mu'$) are full.
It follows that $\mu$ and $\mu'$ are disjoint.

Since the cycles induced by $\mu$ and $\mu'$ are full, it is easy to see that these paths are disconnected.
For this we can argue exactly as we did above for $\mu_p$ and $\mu'_p$.
Thus, by Lemma \ref{commutativity} both paths are non-zero.

This means that in case $\rho_1 = 0$, there is only one term in our linear combination, and we are
done.

Now assume all $\rho_i$ are non-zero. If $n=2$, we are done, so assume $n>2$. 
Then there are at least three terms and 
we have a non-zero scalar $\lambda$ such that $\rho_1 + \lambda \rho_2 = 0$. Hence, we can write
$$k_1 \rho_1 + \cdots k_n \rho_n = 
k_1(\rho_1 + \lambda \rho_2) + (-k_1 \lambda + k_2)\rho_2 + k_3 \rho_3 + \cdots + k_n \rho_n. $$

Since $\rho$ is not in $\r I + I \r$, we must have that either
$\rho_1 + \lambda \rho_2$ or $(-k_1 \lambda + k_2)\rho_2 + k_3 \rho_3 + \cdots + k_n \rho_n$
is minimal. If $\rho_1 + \lambda \rho_2$ is minimal we have an arrow $i \to j$ by Lemmas 
\ref{two-paths} and \ref{only-disc}.
Hence $\rho = k_1 \rho_1 + k_2 \rho_2$, by Lemma \ref{maxtwo}, and we are done.
In case $\rho_1 + \lambda \rho_2$ is not minimal, we have an expression of
$\rho$ as a minimal relation as a sum of a minimal relation with $n-1$ terms plus something
not minimal, and we are done by induction.
\end{proof}

The first part of the next lemma is a general and well-known fact. We include 
a proof for the convenience of the reader.

\begin{lem}\label{generates}
\begin{itemize}
\item[(a)]{The ideal $I$ is generated by minimal relations.}
\item[(b)]{The ideal $I$ is generated by minimal zero-relations and minimal commutativity-relations.}
\end{itemize}
\end{lem}

\begin{proof}
(a): Let $\{\rho_i \}$ be a set of non-zero elements in $I$, and 
let $\{\overline{\rho_i} \}$ be the corresponding set in $I / (\r I+ I \r)$.
We claim that if the set $\{\overline{\rho_i} \}$ generates $I / (\r I+ I \r)$, then the
set $\{\rho_i \}$ generates $I$.
Let $\rho$ be an element of $I$. Then there are scalars $k_1, \dots,k_m$ and 
elements $\rho_1, \dots, \rho_m$ such that 
$\rho = k_1 \rho_1 + \cdots k_m \rho_m + i_1 r_1 + \cdots + i_s r_s + r'_1 i'_1 + \cdots + r'_t i'_t$,
with $i_1, \dots, i_s, i'_1, \dots, i'_t$ in $I$ and
with $r_1, \dots, r_s, r'_1, \dots, r'_t$ in $\r$. Each of the elements 
$i_1, \dots, i_s, i'_1, \dots, i'_t$ can be written as a
linear combination of $\rho_i$'s plus something in $\r I + I \r$. If we iterate this
process $l$ times we will have an expression of $\rho$ as $\rho = \rho' + \rho''$, where
$\rho'$ is in the ideal generated by the $\rho_i$'s while $\rho''$ is in $\r^l$.
If we choose $l$ bigger than the maximum of the lengths of the paths occurring as terms
in $\rho$, then it is clear that the elements in $\rho''$ will cancel with some of the elements in $\rho'$,
and we see that $\rho$ is actually in the ideal generated by the $\rho_i$'s. This finishes the proof of (a).
It is clear that (b) also follows from this proof, using Lemma \ref{sum}.
\end{proof}

\section{Main results}

We now  give the main results of this paper.
Note that by Lemma \ref{alternating}, any arrow lying on an oriented cycle, also lies on a full oriented
cycle.

\begin{thm} 
Let $Q$ be a finite quiver, $K$ an algebraically closed field and $I$ an ideal in the path algebra $KQ$, such that
$\G = KQ/I$ is a cluster-tilted algebra of finite representation type.
\begin{itemize}
\item[(a)]{Let $j$ and $i$ be vertices of $Q$ and assume there is an arrow $i \to j$. Then there are at most two shortest paths from
$j$ to $i$.}
\begin{itemize}
\item[(i)]{If there is exactly one, then this is a minimal zero-relation.}
\item[(ii)]{If there are two, $\rho$ and $\mu$, then $\rho$ and $\mu$ are not zero in $\G$,
they are disconnected, and there is a minimal relation $\rho + \lambda \mu$ for some 
$\lambda \neq 0$ in $K$.}
\end{itemize}
\item[(b)]{Up to multiplication by non-zero elements of $K$, there are no other minimal 
zero-relations or commutativity-relations than those described in (a) for some arrow $i \to j$.}
\item[(c)]{The ideal $I$ is generated by minimal zero-relations and minimal commutativity-relations.}
\end{itemize}
\end{thm}

\begin{proof}
Assume there is an arrow $i \to j$. Consider the paths from $j$ to $i$ which give rise to a full and pure cycle.

Assume first that one of the paths $\rho$ is zero in $\G$.
Since the associated cycle is full and pure, $\rho$ is a minimal zero-relation for the factor algebra
obtained by deleting all vertices not on this cycle, and hence also for $\G$.

Assume now there are at least two distinct 
paths $\rho \colon j= j_0 \to j_1 \to \cdots \to j_n = i$ and 
$\mu \colon j=j'_0 \to j'_1 \to \cdots \to j'_n = i$ from $j$ to $i$, 
each giving rise to a full and pure cycle.
By copying arguments in the proof of Lemma \ref{sum}, we get
that 
$\rho$ and $\mu$ are disconnected, and thus
it follows from Proposition \ref{commutativity} that both paths are non-zero. 

Hence we have seen that if there is more than one path giving rise to a full and pure cycle, all 
such paths must be non-zero. By Lemma \ref{maxtwo} there are then at most two such paths. 
We claim that there cannot be exactly one non-zero path. Assume that there is one non-zero path, then
by Lemma \ref{full-cycles} it is full and the corresponding cycle is pure, 
so by Proposition \ref{counting} it is zero in
some proper factor. Thus there must be a distinct path from $j$ to $i$, 
which is also necessarily full and such that the corresponding cycle is pure, 
by Lemma \ref{full-cycles}. 
If there are exactly two such paths, it follows from Lemma \ref{dimone}
that $\rho + \lambda \mu= 0$ for some $\lambda \neq 0$ in $K$.

Assume there is at least one path from $j$ to $i$.
We have now shown that then there is either exactly one path $\rho$
from $j$ to $i$ inducing a full and pure cycle, and in this case $\rho$ is a minimal 
zero-relation, or there are two paths, $\rho$ and $\mu$, both inducing a full and pure cycle. 
In the latter case we have shown that the two paths are disconnected. Thus by Proposition
\ref{commutativity} we have that $\rho + \lambda \mu$ (with $\lambda \neq 0$ 
in $K$) is a minimal relation. That there are no other minimal zero-relations
follows from Proposition \ref{arrow-back}. That there are no other minimal commutativity
relations follows from Lemma \ref{only-disc}. 
\end{proof}

Up to isomorphism of cluster-tilted algebras, we can by the multiplicative basis theorem \cite{bgrs} assume
that all $\lambda$'s in the above statement are $-1$. We have the following nice
consequence, using Lemma \ref{generates}. Recall that when there is an arrow $i \to j$,
we call a path from $j$ to $i$ shortest if 
it contains no proper subpath
which is a cycle and if the full subquiver generated
by the induced oriented cycle contains no further arrows.

\begin{thm}\label{maintwo}
Let $\G$ be a cluster-tilted algebra of finite representation type 
with quiver $Q$. Let $I'$ be the ideal in $KQ$ with the 
following generators.
For $\alpha \colon i \to j$ choose a generator $\rho$ if there is exactly one shortest path $\rho$ from $j$ to $i$.
Choose a generator $\rho - \mu$ if there are two shortest paths $\rho$ and $\mu$ from $j$ to $i$. 
Then $\G \simeq KQ/I'$.  
\end{thm}

\begin{cor}
A (basic) cluster-tilted algebra of finite representation type is determined by its quiver (up to isomorphism).
\end{cor}

%

We end with an application of our main result to completing the solution of Conjecture 1.1
in \cite{ccs1}, proved in \cite{ccs1} for Dynkin type $A$.

The conjecture is: Let $C= \{u_1, \dots, u_n \}$ be a cluster in a cluster algebra 
of simply-laced Dynkin type and rank $n$. Let $V$ be the set of all cluster variables for this cluster algebra.
Then there is a bijection from the set of indecomposable modules of $A_C$ to $V \setminus C$ given by 
$\alpha \mapsto \omega_{\alpha}$,
such that $\omega_{\alpha} = \frac{P(u_1, \dots, u_n)}{\prod_i u_i^{n_i}}$, where $P$ is a polynomial 
prime to $u_i$ for all $i$ and where $n_i = n_i(\alpha)$ is the multiplicity of the simple module
$\alpha_i$ in the module $\alpha$.

Here the algebra $A_C$ associated with the cluster $C$ is defined in \cite{ccs1} in terms of the associated
quiver $Q_C$. Let $T$ be the tilting object in the cluster category $\C$ corresponding to $C$
by the correspondence of \cite{bmrrt}. Then the quiver of $\End_{\C}(T)^{\op}$ is the same as $Q_C$,
by \cite[Theorem 6.2]{bmr2}. Our Theorem \ref{maintwo} can be formulated to say that $\End_{\C}(T)^{\op}$ is 
isomorphic to $A_C$. A modification of the conjecture, replacing $A_C$ by $\End_{\C}(T)^{\op}$, was
shown in \cite{ccs2}, and independently by Reiten and Todorov, as a consequence of
results in \cite{ccs1} and \cite{bmr1}. Hence the conjecture follows as originally stated.


\begin{thebibliography}{99}

\bibitem[ARS]{ars} Auslander M., Reiten I., Smal\o \ S. O. \emph{Representation Theory of Artin Algebras},
Cambridge studies in advanced mathmatics 36, Cambridge University Press (1995)

\bibitem[BGZ]{bgz} Barot M., Geiss C., Zelevinsky A. \emph{Cluster algebras of finite type and 
symmetrizable matrices},
preprint math.CO/0411341 (2004)

\bibitem[BGRS]{bgrs} Bautista R., Gabriel, P., Roiter A. V., Salmeron, L. \emph{Representation-finite 
algebras and multiplicative bases}, Invent. Math. 81, no. 2, 217--285 (1985)


\bibitem[Bo]{Bong} Bongatz, K. \emph{Algebras and quadratic forms}, J. London Math. Soc. (2) 38, no.3, 461--469 (1983)

\bibitem[BMR1]{bmr1} Buan A., Marsh R., Reiten I. \emph{Cluster-tilted algebras}, preprint math.RT/0402075, 
to appear in Trans. Amer. Math. Soc. (2004)

\bibitem[BMR2]{bmr2} Buan A., Marsh R., Reiten I. \emph{Cluster mutation via quiver representations},
preprint math.RT/0412077, (2004)

\bibitem[BMRRT]{bmrrt} Buan A., Marsh R., Reineke M., Reiten I., Todorov G. 
\emph{Tilting theory and cluster combinatorics}, preprint math.RT/0402054, 2004, 
to appear in Adv. Math. 

\bibitem[CK]{ck} Caldero P., Keller B. 
\emph{From triangulated categories to cluster algebras}, preprint math.RT/0506018 (2005)

\bibitem[CC]{cc} Caldero P., Chapoton F. \emph{Cluster algebras as Hall algebras of quiver representations}, 
preprint math.RT/0410187 (2004)

\bibitem[CCS1]{ccs1} Caldero P., Chapoton F., Schiffler R. \emph{Quivers with relations arising from
clusters ($A_n$ case)}, preprint  math.RT/0401316, to appear in Trans. Amer. Math. Soc. (2004)

\bibitem[CCS2]{ccs2} Caldero P., Chapoton F., Schiffler R. \emph{Quivers with relations and
cluster-tilted algebras}, preprint math.RT/0411238 (2004)


\bibitem[H]{hap} Happel D. \emph{Triangulated categories in the representation theory 
of finite-dimensional algebras}, London Mathematical Society Lecture Note Series, 119. 
Cambridge University Press, Cambridge (1988)

\bibitem[FZ]{fz} Fomin S., Zelevinsky A. \emph{Cluster Algebras I: Foundations}, 
J. Amer. Math. Soc. 15, no. 2, 497--529 (2002)

\bibitem[KR]{kr} Keller B., Reiten I. \emph{Cluster-tilted algebras are Gorenstein and stably Calabi-Yau
}, preprint math.RT/0512471 (2005)

\bibitem[MRZ]{mrz} Marsh R., Reineke M., Zelevinsky A. \emph{Generalized associahedra 
via quiver representations}, Trans. Amer. Math. Soc. 355, no. 1, 4171--4186 (2003)  

\bibitem[R]{r} Ringel C. M., \emph{Tame algebras and integral quadratic forms}, 
Springer Lecture Notes in Mathematics 1099 (1984) 

\bibitem[S]{s} Seven A. \emph{Recognizing cluster algebras of finite type}, preprint math.CO/0406545 (2004)


\end{thebibliography}
\end{document}